**Hybrid of monolithic and staggered solution techniques for the computational analysis of fracture, assessed on fibrous network mechanics**

Vedad Tojaga*, Artem Kulachenko, Sören Östlund, T. Christian Gasser

*Solid Mechanics, Department of Engineering Mechanics, KTH Royal Institute of Technology, SE-100 44 Stockholm, Sweden*

**Abstract**

The computational analysis of fiber network fracture is an emerging field with application to paper, rubber-like materials, hydrogels, soft biological tissue, and composites. Fiber networks are often described as probabilistic structures of interacting one-dimensional elements, such as truss-bars and beams. Failure may then be modeled as strong discontinuities in the displacement field that are directly embedded within the structural finite elements. As for other strain-softening materials, the tangent stiffness matrix can be non-positive definite, which diminishes the robustness of the solution of the coupled (monolithic) two-field problem. Its uncoupling, and thus the use of a staggered solution method where the field variables are solved alternatingly, avoids such difficulties and results in a stable, but sub-optimally converging solution method. In the present work, we evaluate the staggered against the monolithic solution approach and assess their computational performance in the analysis of fiber network failure. We then propose a hybrid solution technique that optimizes the performance and robustness of the computational analysis. It represents a matrix regularization technique that retains a positive definite element stiffness matrix while approaching the tangent stiffness matrix of the monolithic problem. The approach is general and may also accelerate the computational analysis of other failure problems.

## 1     Introduction

The theory of strong discontinuities (e.g., see [1–11]) has made a significant impact on the computational analysis of fracture and extended the functionality of many commercial finite element method (FEM) packages. Strong discontinuities enrich the displacement solution space with jumps across the fracture surfaces. The embedded discontinuity finite element method (ED-FEM or E-FEM) and the extended finite element method (X-FEM) have gained the most popularity (e.g., see [12,13] for comparative studies). The E-FEM and the X-FEM represent strong discontinuities via elemental and nodal enrichments of the displacement solution space, respectively. It allows the approximation of fracture problems using coarse FEM meshes, no re-meshing, and results in mesh-independent simulation results. However, said approaches often require a crack-tracking algorithm (except for recent work [14] on E-FEM without crack tracking) and exhibit difficulties capturing crack branching and coalescence of multiple cracks, difficulties that are avoided with the newly emerged phase-field fracture method (e.g., see [15,16]). Unlike discrete crack approaches, it is a smeared crack approach like damage mechanics [17–20] and therefore requires a localization limiter or a characteristic length scale parameter. Such measure is difficult to determine in practice [21] and can be difficult to implement upon unstructured meshes [22]. In addition, a sufficiently refined mesh is necessary to adequately describe the mechanics of the localization zone, leading to the high computational cost of smeared crack approaches. Thus, it prevents from an efficient analysis of failure in fiber networks where the fracture process zone is to be resolved in the individual fibers.

* Corresponding author. arXiv:2203.15924 **[math.NA]**. CC BY-NC-ND 4.0
Email address: tojaga@kth.se (V. Tojaga) and gasser@kth.se (T.C. Gasser).

The tangent stiffness matrix of the two-field problem can be indefinite or negative definite [23–26]. The energy potential to be minimized is therefore non-convex. Uncoupling the problem and using a staggered solution approach that solves for the displacement and the state of damage alternatingly, overcomes this issue and results in a convex alternating minimization problem. It is the preferred approach in phase-field modeling of fracture [27] and available in commercial software [28–30]. It results in a positive definite tangent stiffness matrix, but does not represent the linearized residuum, and the solution then converges much slower as compared to Newton-Raphson iterations. The efficiency of the consistently linearized coupled (monolithic) problem is therefore sacrificed towards the robustness of the solution technique [31].

The E-FEM has received significant attention in the description of failure in truss-bars, beams, and structures thereof [32–46]. No crack-tracking algorithm is then needed [47]. The E-FEM framework allows static condensation of the additional DOF associated with the fracture kinematics [33] directly at the element level, resulting in an operator-splitting method [34] for the evaluation of these DOF. As with the phase-field method, a staggered solution method may be used, and the continuous and discontinuous fields are minimized alternatively [39]. In the present work, we expand these ideas and propose a hybrid method to optimize the performance and robustness of the E-FEM models in the analysis of fiber network failure [39,48]. The hybrid solution approach is general and may also be applied to other fracture mechanics problems.

## 2    Enhanced finite element formulation

We consider a local Cartesian coordinate system $\{\boldsymbol{e}_x, \boldsymbol{e}_y, \boldsymbol{e}_z\}$ in the description of a 3D Timoshenko beam of the length $L$ and the cross-section $A$. The beam's neutral axis $x \in [0, L]$ is aligned with the $x$-coordinate along which distributed $\boldsymbol{f}$ as well as concentrated $\boldsymbol{F}$ loads are applied. The displacements $\boldsymbol{u}(x)$ and the rotations $\boldsymbol{\theta}(x)$ are collectively represented by the generalized displacement vector $\overline{\boldsymbol{u}} =$ $[u_x \quad u_y \quad u_z \quad \theta_x \quad \theta_y \quad \theta_z]^T$ with the subscripts $x$, $y$ and $z$ denoting the respective displacement and rotation components. The generalized strain measures of the 3D Timoshenko beam then read

$$\left.\begin{aligned}\varepsilon = \frac{du_x}{dx} \; ; \; \gamma_y = \frac{du_y}{dx} - \theta_z \; ; \; \gamma_z = \frac{du_z}{dx} + \theta_y \; ; \\ \kappa_x = \frac{d\theta_x}{dx} \; ; \; \kappa_y = \frac{d\theta_y}{dx} \; ; \; \kappa_z = \frac{d\theta_z}{dx} \; ,\end{aligned}\right\} \quad (1)$$

where $\varepsilon$ is the axial strain, $\gamma_y, \gamma_z$ denote the shear strains along the $y$ and $z$ directions, $\kappa_x$ is the change in the angle of twist $\theta_x$ around the beam's neutral axis and $\kappa_y$, $\kappa_z$ represent bending curvatures. The stress resultants $\boldsymbol{\sigma} = [N \quad Q_y \quad Q_z \quad M_x \quad M_y \quad M_z]^T$ are conjugate to the generalized strain measures $\bar{\boldsymbol{\epsilon}} = [\varepsilon \quad \gamma_y \quad \gamma_z \quad \kappa_x \quad \kappa_y \quad \kappa_z]^T$. Here, the axial force $N$, the shear forces $Q_y, Q_z$ along the $y$ and $z$ directions, the torsional moment $M_x$, and the bending moments $M_y, M_z$ along the $y$ and $z$ directions have been introduced.

The strong form of the local equilibrium equations [49] reads

$$\left.\begin{aligned}\frac{dN}{dx} + q_x = 0 \; ; \; \frac{dQ_y}{dx} + q_y = 0 \; ; \; \frac{dQ_z}{dx} + q_z = 0 \; ; \\ \frac{dM_x}{d_x} + m_x = 0 \; ; \; \frac{dM_y}{dx} - Q_z + m_y = 0 \; ; \; \frac{dM_z}{dx} - Q_y + m_z = 0 \; ,\end{aligned}\right\} \quad (2)$$

where $q_x, q_y, q_z$ and $m_x, m_y, m_z$ denote the components of the distributed forces and moments per unit length, respectively.

For simplicity, we consider a two-node beam element of the length $l_e$ with the linear shape functions $N_1 = 1 - x/l_e$ and $N_2 = x/l_e$, and a single Gauss point in the center is used to integrate the FEM equations. The interpolated generalized displacement field then takes the form

$$\bar{\boldsymbol{u}} = \sum_{a=1}^{2} N_a \bar{\boldsymbol{u}}_a = N_1 \bar{\boldsymbol{u}}_1 + N_2 \bar{\boldsymbol{u}}_2, \quad (3)$$

where $\bar{\boldsymbol{u}}_1$ and $\bar{\boldsymbol{u}}_2$ are the six-dimensional generalized nodal displacement vectors. In addition, a softening hinge allows for the formation of localized failure in the center of the Timoshenko beam. We therefore enrich (3) with a discrete displacement/rotation jump $\boldsymbol{\xi} = \left[\xi_{u_x} \quad \xi_{u_y} \quad \xi_{u_z} \quad \xi_{\theta_x} \quad \xi_{\theta_y} \quad \xi_{\theta_z}\right]^T$ situated in the center $x_c$ of the element. Thus,

$$\boldsymbol{u} = N_1 \bar{\boldsymbol{u}}_1 + N_2 \bar{\boldsymbol{u}}_2 + H_{x_c} \boldsymbol{\xi} \quad (4)$$

represents the total displacements, where

$$H_{x_c}(x) = \begin{cases} 1 & \text{at } x > x_c \\ 0 & \text{at } x \leq x_c \end{cases} \quad (5)$$

denotes the Heaviside step function centered at the middle of the element $x_c$. With the boundary conditions $\boldsymbol{u}_1 = \bar{\boldsymbol{u}}_1$ and $\boldsymbol{u}_2 = \bar{\boldsymbol{u}}_2 + \boldsymbol{\xi}$, where $\boldsymbol{u}_1$ and $\boldsymbol{u}_2$ are the total nodal displacement vectors, the displacement field (4) reads

$$\boldsymbol{u} = \boldsymbol{N}\boldsymbol{d} + (H_{x_c} - N_2)\boldsymbol{\xi}, \tag{6}$$

where $\boldsymbol{d} = [\boldsymbol{u}_1 \quad \boldsymbol{u}_2]^T$ represents the element's nodal displacement vector, and

$$\boldsymbol{N} = \begin{bmatrix} N_1 & 0 & 0 & 0 & 0 & 0 & N_2 & 0 & 0 & 0 & 0 & 0 \\ 0 & N_1 & 0 & 0 & 0 & 0 & 0 & N_2 & 0 & 0 & 0 & 0 \\ 0 & 0 & N_1 & 0 & 0 & 0 & 0 & 0 & N_2 & 0 & 0 & 0 \\ 0 & 0 & 0 & N_1 & 0 & 0 & 0 & 0 & 0 & N_2 & 0 & 0 \\ 0 & 0 & 0 & 0 & N_1 & 0 & 0 & 0 & 0 & 0 & N_2 & 0 \\ 0 & 0 & 0 & 0 & 0 & N_1 & 0 & 0 & 0 & 0 & 0 & N_2 \end{bmatrix} \tag{7}$$

denotes the corresponding interpolation matrix. The strain field then takes the form

$$\boldsymbol{\epsilon} = \underbrace{\boldsymbol{B}\boldsymbol{d} + \boldsymbol{G}\boldsymbol{\xi}}_{\epsilon_{\text{bulk}}} + \delta_{x_c}\boldsymbol{\xi}, \tag{8}$$

and the matrix

$$\boldsymbol{G} = G\boldsymbol{I} \tag{9}$$

interpolates the jump, where $G = -dN_2/dx = -1/l_e$ and $\boldsymbol{I}$ denotes the 6 × 6 identity matrix. In addition, $\delta_{x_c} = dH_{x_c}/dx$ is the Dirac delta function centered at $x_c$, and

$$\boldsymbol{B} = \begin{bmatrix} B_1 & 0 & 0 & 0 & 0 & 0 & B_2 & 0 & 0 & 0 & 0 & 0 \\ 0 & B_1 & 0 & 0 & 0 & -N_1 & 0 & B_2 & 0 & 0 & 0 & -N_2 \\ 0 & 0 & B_1 & 0 & N_1 & 0 & 0 & 0 & B_2 & 0 & N_2 & 0 \\ 0 & 0 & 0 & B_1 & 0 & 0 & 0 & 0 & 0 & B_2 & 0 & 0 \\ 0 & 0 & 0 & 0 & B_1 & 0 & 0 & 0 & 0 & 0 & B_2 & 0 \\ 0 & 0 & 0 & 0 & 0 & B_1 & 0 & 0 & 0 & 0 & 0 & B_2 \end{bmatrix} \tag{10}$$

represents the standard 6 × 12 strain-displacement interpolation matrix, where $B_1 = dN_1/dx = -1/l_e$ and $B_2 = dN_2/dx = 1/l_e$.

Given the admissible variation $\delta\boldsymbol{d}$ of the element nodal displacements, $\delta\boldsymbol{u}$ of the displacement field and $\delta\boldsymbol{\xi}$ of the corresponding jump, the internal and respective external virtual work of the beam formulation reads

$$\left.\begin{aligned} \delta w_{\text{int}} &= \int_{l_e} \delta\boldsymbol{\epsilon}^T \boldsymbol{\sigma}\, dx = \delta\boldsymbol{d}^T \int_{l_e} \boldsymbol{B}^T \boldsymbol{\sigma}\, dx + \delta\boldsymbol{\xi}^T \left(\int_{l_e} \boldsymbol{G}\boldsymbol{\sigma}\, dx + \int_{l_e} \delta_{x_c}\boldsymbol{\sigma}\, dx\right) \,; \\ \delta w_{\text{ext}} &= \delta\boldsymbol{u}^T(x_F)\,\boldsymbol{F} + \int_{l_e} \delta\boldsymbol{u}^T \boldsymbol{f}\, dx = \delta\boldsymbol{d}^T\left(\boldsymbol{N}^T(x_F)\,\boldsymbol{F} + \int_{l_e} \boldsymbol{N}^T \boldsymbol{f}\, dx\right) + \delta\boldsymbol{\xi}^T\bigl(H_{x_c} - N_2(x_F)\bigr)\boldsymbol{F}\,, \end{aligned}\right\} \tag{11}$$

where $\boldsymbol{f} = [q_x\ q_y\ q_z\ m_x\ m_y\ m_z]^T$ represents a distributed load along the beam and $\boldsymbol{F}$ a concentrated load applied at the position $x_F$ along the beam. In the derivation of $(11)_2$, we used the identity $\int_{l_e}(H_{x_c} - N_2)\,dx = 0$, a consequence of having the discontinuity in the middle of the beam, $x_c = l_e/2$. In addition, the second term of $(11)_2$ vanishes if $x_F = 0$ or $x_F = l_e$, which we apply here.

From the principle of virtual work, $\delta w_{\text{int}} - \delta w_{\text{ext}} = 0$, and the independence of the admissible variations $\delta \boldsymbol{d}$, $\delta \boldsymbol{\xi}$, we obtain the two variational statements [33,34,50]

$$\left.\begin{array}{r}\delta \boldsymbol{d}^T \left( \int_{l_e} \boldsymbol{B}^T \boldsymbol{\sigma}\, dx - \boldsymbol{N}^T(x_F)\boldsymbol{F} + \int_{l_e} \boldsymbol{N}^T \boldsymbol{f}\, dx \right) = 0 \ ; \\ \delta \boldsymbol{\xi}^T \left( \int_{l_e} \boldsymbol{G} \boldsymbol{\sigma}\, dx + \boldsymbol{t} \right) = 0 \ , \end{array}\right\} \tag{12}$$

where the condition $\boldsymbol{t} = \int_{l_e} \delta_{x_c} \boldsymbol{\sigma}\, dx$ has been used. The condition $(12)_2$ represents the local residual at the discontinuity, $x = x_c$, and ensures the equilibrium between the traction $\boldsymbol{t}$ acting on the discontinuity and the stress resultant vector $\boldsymbol{\sigma}$ in the bulk material, $-\boldsymbol{\sigma} + \boldsymbol{t} = \boldsymbol{0}$.

The implementation of (12) results in the set

$$\boldsymbol{r}_d^e = \boldsymbol{f}^{e\,\text{int}} - \boldsymbol{f}^{e\,\text{ext}} = \boldsymbol{0} \ ; \quad \boldsymbol{r}_\xi^e = \int_{l_e} \boldsymbol{G}\boldsymbol{\sigma}\, dx + \boldsymbol{t} = \boldsymbol{0} \tag{13}$$

of non-linear equations at the element level, where $\boldsymbol{f}^{e\,\text{int}} = \int_{l_e} \boldsymbol{B}^T \boldsymbol{\sigma}\, dx$ and $\boldsymbol{f}^{e\,\text{ext}} = \int_{l_e} \boldsymbol{N}^T \boldsymbol{f}\, dx + \boldsymbol{N}^T(x_F)\boldsymbol{F}$ are the internal and external element nodal force vectors, respectively.

Incremental formulations are used to express the constitutive relations of the bulk material and the fracture process zone, expressions that close the set of equations (13). We use

$$\Delta \boldsymbol{\sigma} = \boldsymbol{C} \Delta \boldsymbol{\epsilon}_{\text{bulk}} \ ; \quad \Delta \boldsymbol{t} = \boldsymbol{H} \Delta \boldsymbol{\xi} \tag{14}$$

to describe the development of the stress resultant $\boldsymbol{\sigma}$ and the traction $\boldsymbol{t}$, where $\boldsymbol{C}$ and $\boldsymbol{H}$ denote the respective tangent constitutive tensors. Given a beam of the cross-section $A$ that is made of a linear-elastic material with Young's modulus $E$ and shear modulus $G$, the tangent

$$\boldsymbol{C} = \begin{bmatrix} EA & 0 & 0 & 0 & 0 & 0 \\ 0 & kGA & 0 & 0 & 0 & 0 \\ 0 & 0 & kGA & 0 & 0 & 0 \\ 0 & 0 & 0 & GJ & 0 & 0 \\ 0 & 0 & 0 & 0 & GI_{11} & 0 \\ 0 & 0 & 0 & 0 & 0 & GI_{22} \end{bmatrix} \tag{15}$$

determines the development of $\boldsymbol{\sigma}$, where $k$ is the shear correction factor, $J$ is the polar moment of inertia, whilst $I_{11}$ and $I_{22}$ denote the area moments of inertia. Note that (15) is the simplest set of uncoupled linear elastic constitutive equations and based on the assumption that the beam cross-section possesses appropriate symmetries [51].

With the strain $\epsilon_{\text{bulk}}$ in the bulk material (8), the increment of the stress resultant

$$\Delta \boldsymbol{\sigma} = \boldsymbol{C}(\boldsymbol{B}\Delta \boldsymbol{d} + \boldsymbol{G}\Delta \boldsymbol{\xi}), \tag{16}$$

and incremental local equilibrium across the discontinuity $-\Delta\boldsymbol{\sigma} + \Delta\boldsymbol{t} = \boldsymbol{0}$ are to be enforced.

The linearization of the residual force equations (13) with respect to the unknown displacement $\Delta \boldsymbol{d}$ and the discontinuous displacement $\Delta \boldsymbol{\xi}$ yields the system

$$\begin{bmatrix} \boldsymbol{K}^e_{dd} & \boldsymbol{K}^e_{d\xi} \\ \boldsymbol{K}^e_{\xi d} & \boldsymbol{K}^e_{\xi\xi} \end{bmatrix} \begin{bmatrix} \Delta \boldsymbol{d} \\ \Delta \boldsymbol{\xi} \end{bmatrix} = \begin{bmatrix} \boldsymbol{f}^{e\,\text{int}} - \boldsymbol{f}^{e\,\text{ext}} \\ \boldsymbol{0} \end{bmatrix} \tag{17}$$

with the sub-matrices

$$\left. \begin{aligned} \boldsymbol{K}^e_{dd} &= \int_{l_e} \boldsymbol{B}^T \boldsymbol{C} \boldsymbol{B}\, dx \;\;;\;\; \boldsymbol{K}^e_{d\xi} = \int_{l_e} \boldsymbol{B}^T \boldsymbol{C} \boldsymbol{G}\, dx \;\;;\; \\ \boldsymbol{K}^e_{\xi d} &= \int_{l_e} \boldsymbol{G}\boldsymbol{C}\boldsymbol{B}\, dx + \boldsymbol{C}^*\boldsymbol{B} \;\;;\;\; \boldsymbol{K}^e_{\xi\xi} = \int_{l_e} \boldsymbol{G}^T \boldsymbol{C} \boldsymbol{G}\, dx + \boldsymbol{H} \end{aligned} \right\} \tag{18}$$

The term $\boldsymbol{C}^*\boldsymbol{B}$ is explained in **Section 3**.

In the implementation of these equations, we distinguish between elastic and failure loading. We follow the framework of inelasticity [23,52] and the decision is based on the introduction of a failure surface $\Phi$, see Section 3. Given elastic loading, or unloading, the increment $\Delta \boldsymbol{\xi} = \boldsymbol{0}$ and the system reduces to the standard Timoshenko beam FEM model $\boldsymbol{K}^e_{dd}\Delta\boldsymbol{d} = \boldsymbol{f}^{e\,\text{int}} - \boldsymbol{f}^{e\,\text{ext}}$. Different implementations concerning damage loading are discussed in the forthcoming sections.

## 2 Monolithic, staggered and hybrid FEM implementation

The monolithic implementation considers the consistently derived finite element stiffness as shown in (17), where off-diagonal terms couple the increment of the displacement $\Delta \boldsymbol{d}$ and the increment of the jump $\Delta \boldsymbol{\xi}$. As with other embedded approaches, our model allows for the static condensation of $\Delta \boldsymbol{\xi}$ directly at the element level, which then results in a displacement-based FEM implementation. The second equation in (17), the internal equilibrium, is then used to substitute $\Delta \boldsymbol{\xi}$ through

$$\Delta \boldsymbol{\xi} = -\left(\boldsymbol{K}^e_{\xi\xi}\right)^{-1} \boldsymbol{K}^e_{\xi d} \Delta \boldsymbol{d}. \tag{19}$$

In the limit $l_e \to 0$, $\boldsymbol{K}^e_{\xi\xi}$ remains positive definite and invertible [6,33]. In general, however, the condition (19) poses a limit on the size $l_e$ of the finite element in the analysis of strain-softening materials. The substitution of $\Delta \boldsymbol{\xi}$ by (19) in (17) results in the system of equations

$$\boldsymbol{K}^e_{\text{mono}}\Delta\boldsymbol{d} = \boldsymbol{f}^{e\text{ int}} - \boldsymbol{f}^{e\text{ ext}} \tag{20}$$

at the element level, where

$$\boldsymbol{K}^e_{\text{mono}} = \boldsymbol{K}^e_{dd} - \boldsymbol{K}^e_{d\xi}\left(\boldsymbol{K}^e_{\xi\xi}\right)^{-1}\boldsymbol{K}^e_{\xi d} \tag{21}$$

denotes the element stiffness. We emerged at a displacement-based model that may now be implemented through the standard user element interface of FEM packages. Whilst the aforementioned solution uses the consistent linearization of the residual forces, it results in a non-positive definite finite element stiffness $\boldsymbol{K}^e_{\text{mono}}$, which then materializes through poor robustness of the monolithic implementation.

Towards reinforcing the robustness of the model, we may uncouple $\boldsymbol{d}, \boldsymbol{\xi}$ and solve the problem in a staggered way. The internal equilibrium equation $(13)_2$ is then solved at the nodal displacement $\boldsymbol{d}_n$ from the previous solution. It explicitly reads

$$\boldsymbol{r}^e_\xi = \int_{l_e} \boldsymbol{G}\boldsymbol{\sigma}(\boldsymbol{d}_n)\,dx + \boldsymbol{t} = \boldsymbol{0}, \tag{22}$$

and the embedded formulation again allows to solve it directly at the element level. The system

$$\boldsymbol{K}^e_{\text{stagg}}\Delta\boldsymbol{d} = \boldsymbol{f}^{e\text{ int}} - \boldsymbol{f}^{e\text{ ext}} \tag{23}$$

of FEM equations may then be assembled, and the solution of the global system yields the nodal displacements. Here,

$$\boldsymbol{K}^e_{\text{stagg}} = \boldsymbol{K}^e_{dd} \tag{24}$$

determines the corresponding finite element stiffness, and given it represents an inconsistently linearized residuum, which results in poor convergence of the staggered approach.

Towards optimizing performance and robustness of the finite element model, we propose the hybrid definition

$$\boldsymbol{K}^e_{\text{hyb}} = \beta \boldsymbol{K}^e_{\text{mono}} + (1-\beta)\boldsymbol{K}^e_{\text{stagg}} \tag{25}$$

of the finite element stiffness, where $\beta$ is a numerical parameter, chosen to ensure a positive definite finite element stiffness matrix $\boldsymbol{K}^e_{\text{hyb}}$. It is set according to

$$0 \leq \beta < \beta_{\text{critical}} \quad \text{such that} \quad \det \boldsymbol{K}^e_{\text{hyb}} > 0, \tag{26}$$

where the condition for $\beta_{\text{critical}} \Rightarrow \det \boldsymbol{K}^e_{\text{hyb}} = 0$ is derived as follows. Excluding the six rigid body motion-related DOFs from the system, the effective stiffness matrices are of the dimension $6 \times 6$, and the only physical root of $\det \boldsymbol{K}^e_{\text{hyb}} = 0$ results in an expression for $\beta_{\text{critical}}$ and (26) then guarantees a positive definite finite element stiffness matrix. Although our stiffness matrices are sparsely populated, the direct solution of $\det \boldsymbol{K}^e_{\text{hyb}} = 0$ requires the eigenvalue analysis of one $12 \times 12$ matrix at every Gauss point for each solution step; a faster and tailored implementation for beam rupture is discussed in **Section 4**.

## 3    Predictor-corrector implementation of beam rupture

Aiming at modeling the failure of soft fibers, we limit ourselves to the description of failure under tension. Therefore, only the component $\xi_{u_x} = \xi$ of the jump $\boldsymbol{\xi}$ is allowed to evolve, whilst $\xi_{u_y} = \xi_{u_z} = \xi_{\theta_x} = \xi_{\theta_y} = \xi_{\theta_z} = 0$. The development of the traction according to $(14)_2$ is then determined by the tangent constitutive tensor component $H_{11} = H < 0$, whilst all the other components of $\boldsymbol{H}$ are 0. As a result, $\boldsymbol{C}^*$ in (18) takes the form

$$\boldsymbol{C}^* = \begin{bmatrix} 0 & 0 & 0 & 0 & 0 & 0 \\ 0 & kGA & 0 & 0 & 0 & 0 \\ 0 & 0 & kGA & 0 & 0 & 0 \\ 0 & 0 & 0 & GJ & 0 & 0 \\ 0 & 0 & 0 & 0 & GI_{11} & 0 \\ 0 & 0 & 0 & 0 & 0 & GI_{22} \end{bmatrix} \quad (27)$$

because we only allow failure in tension. In (16), $\boldsymbol{C}^*\boldsymbol{B}$ is the tangent when $\Delta\xi_{u_y} = \Delta\xi_{u_z} = \Delta\xi_{\theta_x} = \Delta\xi_{\theta_y} = \Delta\xi_{\theta_z} = 0$.

The failure criterion [34]

$$\Phi = N - (\overline{N} + H\alpha) \leq 0 \quad (28)$$

in the stress resultant space, determines the beam's loading condition. Here, $\overline{N}$ is the elastic limit resultant (force), whilst $H < 0$ and $\alpha > 0$ denote the softening modulus and an internal softening variable, respectively. Elastic deformation of the beam is then characterized by $\Phi < 0$, and a loading state that reaches the failure surface, $\Phi = 0$, results in the accumulation of failure. At complete rupture, $N = 0$, the linear softening law (28) yields $\alpha_{\max} = \overline{N}/|H|$ and determines the fracture energy of $G_f = \alpha_{\max}|H|\alpha_{\max}/2 = \overline{N}^2/(2|H|)$.

Towards closing the failure description, the evolution of the internal softening variable is to be linked to the evolution of the jump displacement, and

$$\begin{array}{ll} \Delta\alpha = 0 & \text{at} \quad \phi < 0 \;; \\ \Delta\alpha = \Delta\xi & \text{at} \quad \phi = 0 \end{array} \quad (29)$$

specifies said correspondence at the cases of elastic and failure loading, respectively.

With (16), the internal equilibrium $-\boldsymbol{\sigma} + \boldsymbol{t} = \boldsymbol{0}$ for the only non-trivial loading mode reads $-N + t = 0$, resulting in

$$N = \bar{N} + H\xi = EA(\varepsilon + G\xi). \tag{30}$$

Here, $\epsilon_{\text{bulk}} = \varepsilon + G\xi$ denotes the strain in the bulk material, and $\xi$ is the jump displacement with the initial condition $\xi = 0$ at $N = \bar{N}$. Given $\xi_n$ and $\alpha_n$ from the previous time point, we can iteratively derive $\xi$ and $\alpha$ at the current time point. We therefore expand (30) towards

$$N = EA(\varepsilon + G\xi) = N^{\text{trial}} + \underbrace{EAG\Delta\gamma\,\text{sign}(N)}_{N^{\text{corr}}}, \tag{31}$$

where $N^{\text{trial}} = EA(\varepsilon + G\xi_n)$ and $\Delta\gamma = \xi - \xi_n$ denotes a consistency parameter that ensures (28); it eventually corrects the trial state by $N^{\text{corr}}$. The implementation follows the classical concept of computational inelasticity [53], and once $\Delta\gamma$ is given, $\xi = \xi_n + \Delta\gamma\,\text{sign}(N)$ updates the solution.

Alternatively, equation (31) may be expressed as $(|N| - EAG\Delta\gamma)\,\text{sign}(N) = |N^{\text{trial}}|\,\text{sign}(N^{\text{trial}})$, and with $\Delta\gamma > 0$, $G < 0$ and $EA > 0$, the term in the bracket to the left is always positive. Therefore, $\text{sign}(N) = \text{sign}(N^{\text{trial}})$ and $|N| = |N^{\text{trial}}| + EAG\Delta\gamma$ follows.

With (31), the failure surface (28) then reads

$$\Phi = N - (\bar{N} + H\alpha) = \Phi^{\text{trial}} + (EAG\,\text{sign}(N) - H)\Delta\gamma, \tag{32}$$

where $\Phi^{\text{trial}} = N^{\text{trial}} - (\bar{N} + H\alpha_n)$ describes the failure surface assuming an elastic load increment. At failure loading, $\Phi = 0$, the algorithmic consistency parameter therefore reads

$$\Delta\gamma = \frac{\Phi^{\text{trial}}}{H - EAG\,\text{sign}(N)} > 0 \tag{33}$$

and allows us to update the solution. The algorithmic consistency parameter $\Delta\gamma$ must be positive, a standard stability condition to be enforced to avoid snap-back behavior. **Table 1** summarizes the predictor-corrector implementation.

**Table 1**: Predictor-corrector implementation of beam rupture.

---

1. Compute elastic trial force and test for failure loading

$$N^{\text{trial}} = EA(\varepsilon + G\xi_n)$$

$$\Phi^{\text{trial}} = N^{\text{trial}} - (\bar{N} + H\alpha_n)$$

IF $\Phi^{\text{trial}} \leq 0$ THEN

$$N = N^{\text{trial}} \; ; \; \xi = \xi_n \; ; \; \alpha = \alpha_n$$

ELSE

2. Return mapping

DO WHILE $\Phi^{\text{trial}} >$ THRESHOLD

$$\Delta\gamma = \frac{\Phi^{\text{trial}}}{H - EAG \, \text{sign}(N)}$$

$$N = N^{\text{trial}} + EAG\Delta\gamma \, \text{sign}(N^{\text{trial}})$$

$$\xi = \xi_n + \Delta\gamma \, \text{sign}(N^{\text{trial}})$$

$$\alpha = \alpha_n + \Delta\gamma$$

END DO

END IF

---

## 4 Implementation

The description of beam rupture through failure under tension and the uncoupled constitutive model (15) result in an uncoupling of the condition (25). Consequently, the identification of the stability parameter (26) yields the single scalar equation

$$\beta = \begin{cases} \beta = 1 & \text{if} \quad K^e_{\text{mono } 11} > K_{\text{min}} \; ; \\ \dfrac{K_{\text{min}} - K^e_{\text{stagg } 11}}{K^e_{\text{mono } 11} - K^e_{\text{stagg } 11}} & \text{if} \quad K^e_{\text{mono } 11} < K_{\text{min}} \; , \end{cases} \qquad (34)$$

and avoids then the eigenvalue analysis of the stiffness matrices $\boldsymbol{K}^e_{\text{mono}}$ and $\boldsymbol{K}^e_{\text{stagg}}$. Here, $K_{\text{min}} = h_{\text{tol}}EA/l_e > 0$ determines a minimum stiffness, where the numerical tolerance level $h_{\text{tol}}$ is determined by the precision level of the hardware/software realization, whilst $l_e$ is the element length. To avoid ill-conditioning towards the development of complete fiber rupture, the corresponding elements in the global stiffness matrix are limited to be larger than the minimum stiffness $K_{\text{min}}$, where again said stiffness threshold has been used.

Given the finite element stiffness matrix (25) for the cases of a monolithic ($\beta = 1$), staggered ($\beta = 0$), or hybrid ($0 < \beta < 1$) implementation, we may now rotate the DOF vector $\boldsymbol{d}$, the nodal force vectors $\boldsymbol{f}^{e\,\text{int}}$, and the stiffens $\boldsymbol{K}^e$ into the global Cartesian coordinate system $\{\boldsymbol{e}'_x, \boldsymbol{e}'_y, \boldsymbol{e}'_z\}$ [54]

$$\boldsymbol{d} = \boldsymbol{T}\boldsymbol{d}' \quad ; \quad \boldsymbol{f}^{e\,\text{int}'} = \boldsymbol{T}^T \boldsymbol{f}^{e\,\text{int}} \quad ; \quad \boldsymbol{K}^{e'} = \boldsymbol{T}^T \boldsymbol{K}^e \boldsymbol{T}. \tag{35}$$

Here,

$$\boldsymbol{T} = \begin{bmatrix} \Lambda & 0 & 0 & 0 \\ 0 & \Lambda & 0 & 0 \\ 0 & 0 & \Lambda & 0 \\ 0 & 0 & 0 & \Lambda \end{bmatrix} \tag{36}$$

denotes the corresponding transformation matrix, where $\Lambda_{ij} = \boldsymbol{e}_i^T \boldsymbol{e}'_j; i,j = x,y,z$ are the directional cosines between the global $\{\boldsymbol{e}'_x, \boldsymbol{e}'_y, \boldsymbol{e}'_z\}$ and the local $\{\boldsymbol{e}_x, \boldsymbol{e}_y, \boldsymbol{e}_z\}$ systems, respectively.

**Table 2** summarizes the implementation in the commercial software package ANSYS through the user-element interface.

1. Load $\alpha_n, \xi_n$
2. Global-to-local transformation $\boldsymbol{d} = \boldsymbol{T}\boldsymbol{d}'$
3. Compute $\boldsymbol{\varepsilon} = \boldsymbol{B}\boldsymbol{d}$
4. Compute $Q_y, Q_z, M_x, M_y, M_z$ as $\boldsymbol{\sigma} = \boldsymbol{C}\boldsymbol{\varepsilon}$
5. Element degrees of freedom (DOF) deletion
    IF $(\alpha_n > \alpha_{max})$ THEN
        $K^e_{1j} = K^e_{j1} = K^e_{7j} = K^e_{j7} = 0 \quad ; \quad j = 1, \dots, 12 \quad ; \quad N = 0$
        GOTO 8
    ENDIF
6. Identify $\alpha, \xi, N$ according to **Table 1**
    IF $(\alpha_n > \alpha_{max})$ THEN
        $K^e_{1j} = K^e_{j1} = K^e_{7j} = K^e_{j7} = 0 \quad ; \quad j = 1, \dots, 12 \quad ; \quad N = 0$
        GOTO 8
    ENDIF
7. IF (Monolithic) THEN $\beta = 1$
    ELSEIF (Staggered) THEN $\beta = 0$
    ELSE (Hybrid) THEN $K_{\min} = \frac{h_{\text{tol}}EA}{l_e} \quad ; \quad \beta = \frac{K_{\min} - K^e_{\text{stagg 11}}}{K^e_{\text{mono 11}} - K^e_{\text{stagg 11}}}$
    ENDIF
    Compute element stiffness $\boldsymbol{K}^e = \beta \boldsymbol{K}^e_{\text{mono}} + (1-\beta) \boldsymbol{K}^e_{\text{stagg}}$
8. Compute nodal force vector $\boldsymbol{f}^{e\,\text{int}} = \int_{l_e} \boldsymbol{B}^T \boldsymbol{\sigma} \, dx$
9. Transform to global coordinates $\boldsymbol{K}^{e'} = \boldsymbol{T}^T \boldsymbol{K}^e \boldsymbol{T}$ and $\boldsymbol{f}^{e\,\text{int}'} = \boldsymbol{T}^T \boldsymbol{f}^{e\,\text{int}}$
10. Return $\alpha, \xi, \boldsymbol{K}^{e'}, \boldsymbol{f}^{e\,\text{int}'}$

## 5 Benchmark exercises

### 5.1 Mesh independence

To illustrate the behavior of the finite element, we consider a $0.1\ mm$ long cantilever beam that is loaded in tension. Young's modulus $E = 1.0\ MPa$, the cross-section area $A = 1.0\ mm^2$, and the ultimate tensile force $\bar{N} = 1.0\ N$ further describe the problem. The cantilever beam is discretized with one and ten evenly spaced finite elements, where the ultimate tensile force is lowered by 1% in the most left element. Otherwise, the homogenous stress field of the problem would have resulted in an ambiguous failure pattern, see Section 3 elsewhere [39]. **Figure 1** presents the force-displacement response of our structural problem for three different fracture energies $G_f$, a property also expressed as the area under the curve. As expected from the fracture model, the force-displacement response is independent of the finite element discretization; it is identical between the one-element and ten-element discretization.

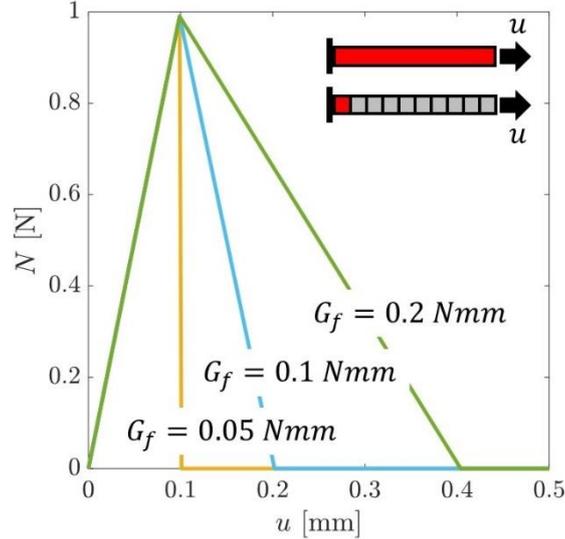

**Fig. 1**: Failure of a cantilever beam under tension using a one-element and ten-element discretization. Axial force $N$ as a function of the axial displacement $u$ at the free end is shown for three different fracture energies $G_f$. The inset illustrates the investigated discretization with the strain localization in the red element.

### 5.2 Tensile test of a fibrous tissue specimen

A planar and random network of interconnected 3D beams describes the mechanics of a fiber network of the densities $\rho_s = \{300; 500; 1000\}\ kg/m^3$ and with the fiber properties listed in **Table 3**. In-plane it covers the area $W \times H = 18 \times 6\ mm^2$ and the displacement $\delta_0 = 9.0\ mm$ is prescribed along one edge, while all six DOFs at the opposite edge are fixed. Out-of-plane displacements and rotations are prevented, and all information on how the mesh has been generated is reported elsewhere [55] together with the ANSYS input file [39]. A quasi-static failure analysis was computed through the incremental application of the prescribed displacement $\delta_0$.

**Figure 2** shows the stress-strain response of the fiber networks. Computational results refer to the staggered solutions. As expected, we observe a more dissipative response of the fiber networks with the higher fracture energy of the fibers. As an example, **Figure 3** shows the evolution of the jump $\xi$ for the densest network of fibers with the fracture energy $G_f = 0.1\ Nmm$. The fiber network configurations refer to the points A, B and C in **Figure 2**. Given said fracture energy, the linear softening law (28) provides the jump $\xi_{max} = 0.85\ mm$ at the state of complete rupture.

**Table 2**: Fiber properties estimated from [56].

| Fiber properties | |
| --- | --- |
| Young's modulus $E$ | $6500\ MPa$ |
| Shear modulus $G$ | $3250\ MPa$ |
| Shear correction factor $k$ | 0.84 |
| Length × Width × Heigth | $2.5 \times \sqrt{0.00028} \times \sqrt{0.00028}\ mm$ |
| Breaking force $\bar{N}$ | $0.2352\ N$ |
| Fracture energy $G_f$ | $0.1\ Nmm$ or $0.2\ Nmm$ |
| Fiber density $\rho_f$ | $1500\ kg/m^3$ |

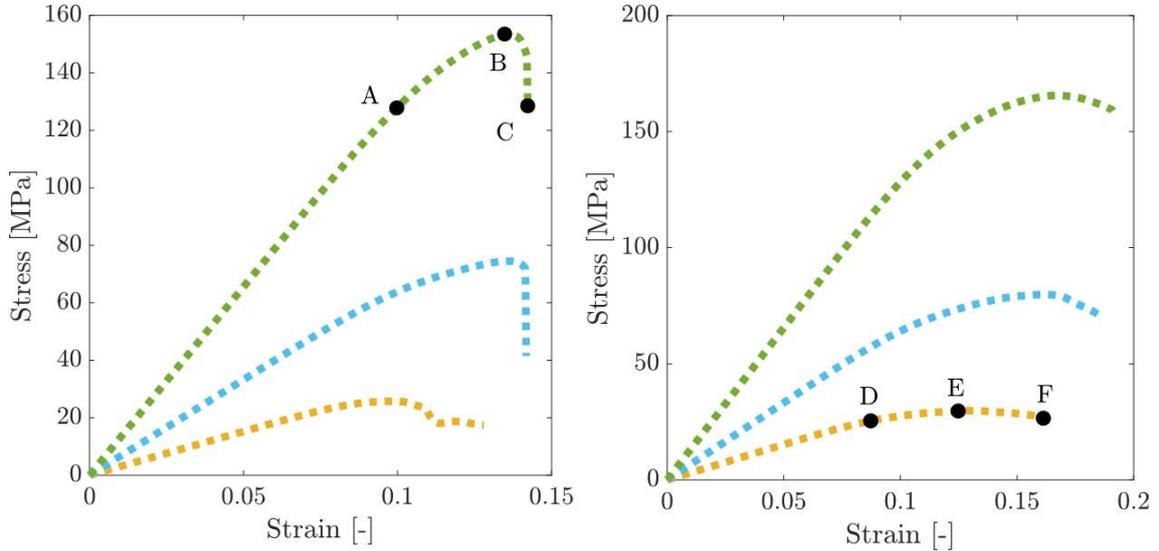

**Fig. 2**: Stress-strain response of fiber networks with the density of $300\ kg/m^3$ (yellow), $500\ kg/m^3$ (blue) and $1000\ kg/m^3$ (green). Dotted lines represent results achieved with the staggered solution method. The fibers within the network are modeled as Timoshenko beams with the properties listed in **Table 3** and the facture energy $G_f = 0.1\ Nmm$ (left column) and $G_f = 0.2\ Nmm$ (right column).

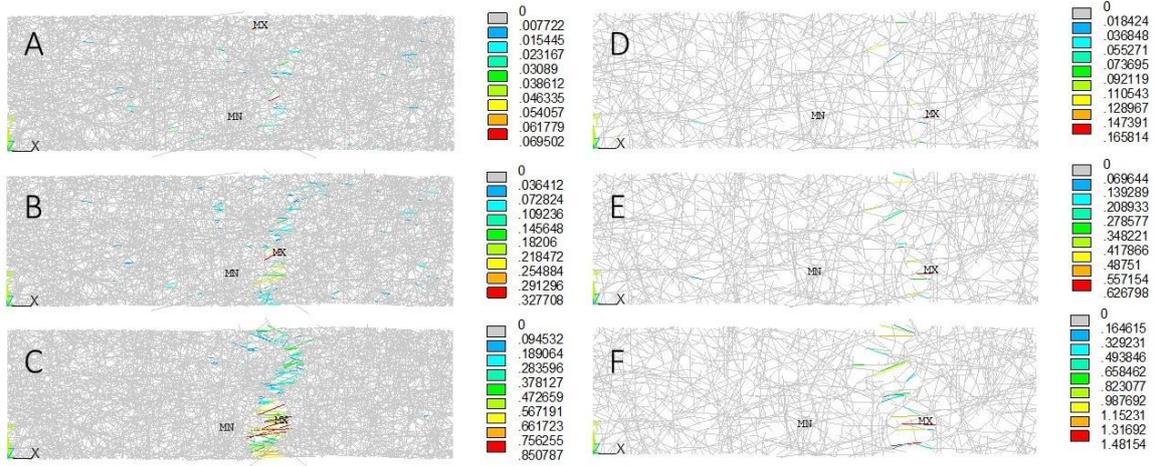

**Fig. 3**: Illustration of the development of failure in the fibrous network of the density $1000\ kg/m^3$ (left column) and $300\ kg/m^3$ (right column). The fibers within the network are modeled as Timoshenko beams with the properties listed in **Table 3** and the facture energy $G_f = 0.1\ Nmm$ (left column) and $G_f = 0.2\ Nmm$ (right column). The jump $\xi\ [mm]$ in the displacement field is shown at points A, B, and C (left column) and D, E, and F (right column), highlighted in **Fig. 2**. At complete rupture, the jump $\xi_{max} = 0.85\ mm$ (left column) and $\xi_{max} = 1.70\ mm$ (right column). Any $\xi > \xi_{max}$ indicates the fiber has failed.

### 5.2.1 Numerical stability

The numerical stability of the $1000\ kg/m^3$ dense fibrous network made of fibers with the fracture energy of $0.2\ Nmm$ was investigated with respect to the size of the displacement increment. The prescribed displacement was applied through 500, 2k, 4k, and 8k steps. In each sub-step a maximum of 500 equilibrium iterations were allowed (NEQIT=500 in ANSYS), before the next step was processed. Whilst the staggered solution converged for all step sizes, the monolithic failed for some step sizes before the staggered solution. **Figure 4** reports the results from the stability analysis and indicates the points of failure of the monolithic solution. Failure is most likely linked to the inability to minimize the related non-convex problem, and a further investigation of ANSYS' internal algorithm was not feasible. Given the monolithic method did not fail, the monolithic and staggered approaches solve the same equilibrium equations (13) and therefore result in the same stress-stain response of the fibrous network.

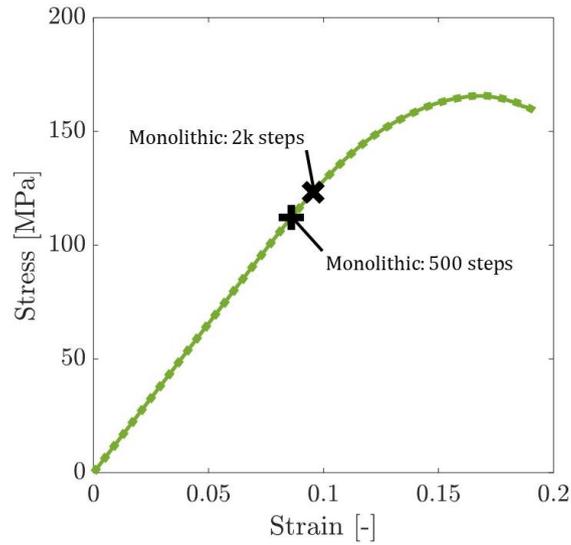

**Fig. 4**: Numerical stability with respect to the size of the displacement increment in the failure analysis of the 1000 $kg/m^3$ dense fibrous network. The fibers within the network are modeled as Timoshenko beams with the properties listed in **Table 3** and the facture energy $G_f = 0.2\ Nm$. Monolithic (solid line) and staggered (dotted line) solutions are explored, where a prescribed displacement was applied through 500, 2k, 4k, and 8k steps, respectively. Crosses denote the point of termination of the (monolithic) computation.

### 5.3 Tensile test of a notched fibrous tissue specimen

The afore explored specimen geometry (see **Section 5.2**) is modified, and a sharp notch with an opening angle of 20° is introduced in the center of the tensile specimen. A 1000 $kg/m^3$ dense fibrous network with the fiber properties listed in **Table 3** is considered. **Figure 5** shows the stress-strain response of the notched tensile specimen together with the evolution of the jump $\xi$. Aligned with the previous problem, the monolithic approach is numerically unstable, given the prescribed displacement was applied through 2k and 10k steps, respectively.

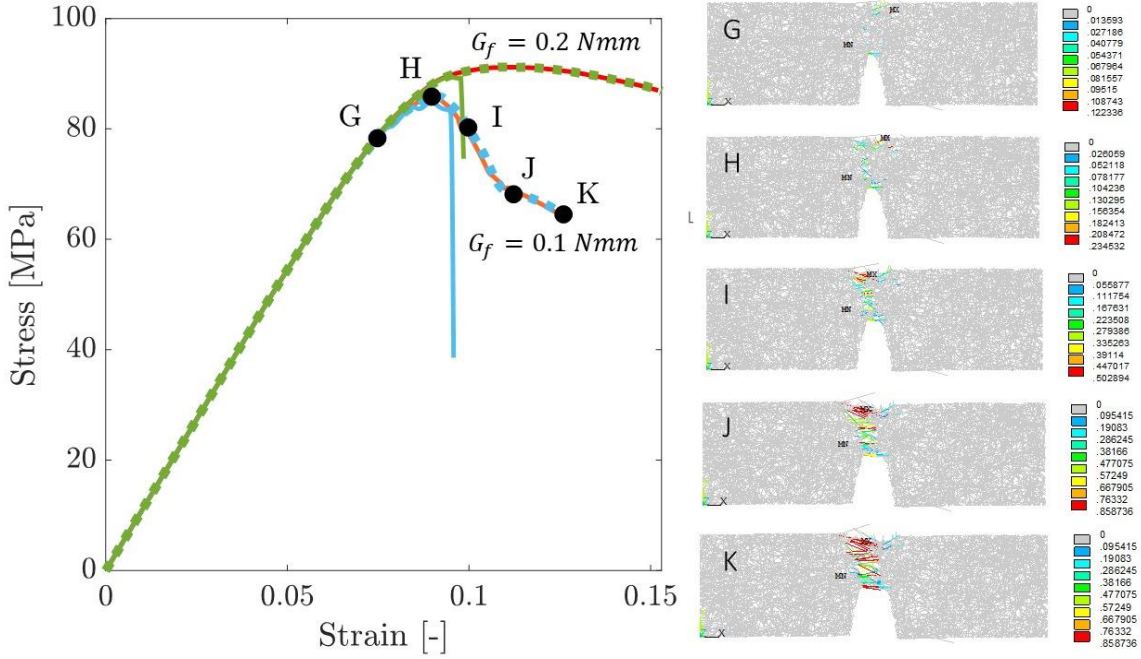

**Fig. 5**: Left: Stress-strain response of fiber network with the density of $1000 \, kg/m^3$ and fibers modeled as Timoshenko beams with the properties listed in **Table 3**. Solid and dotted lines represent results achieved with the monolithic and the staggered solution method, with the monolithic approach failing for large displacement steps. For illustration purposes, the green-red and the blue-orange transition of the solid lines denote a decrease in the time step size for the monolithic solutions to converge to the staggered solutions. Right: Development of failure in the fibrous network with the jump $\xi_{max} = 0.85 \, mm$ at complete rupture. Any $\xi > \xi_{max}$ indicates the fiber has failed.

## 6    Performance and stability of the hybrid solution technique

Among all our simulations, the notched fibrous tissue specimen made of beams with the fracture energy $G_f = 0.1 \, Nmm$ was computationally most demanding. This case will therefore be considered to explore the performance of the proposed hybrid solution technique (25) and benchmarking it against the monolithic and staggered schemes. Towards the optimization of the computation time of large fiber networks with many DOFs, we limit the maximum number of steps to 500 in our benchmarking exercise. We emphasize that at small displacement increments even the staggered solution method converges within a low number of iterations (2 to 4) and then no gain in performance is to be achieved with the hybrid solution. **Table 4** lists the cumulative iterations, the number of all iterations needed to compute the solution until the prescribed displacement $\delta_0 = 1.08 \, mm$ (corresponding to the endpoint in **Fig. 5**) is reached, a measure sensitive to the computational effort to solve the problem. In addition, **Figure 8** shows the resulting stress-strain curves and the predicted end-configurations for the hybrid solution technique, similar to the results shown in **Fig. 5**.

**Table 4**: Cumulative iterations (ANSYS variable CUM ITER) to analyze the failure of the notched fibrous tissue specimen formed by beams with the fracture energy $G_f = 0.1\ Nmm$. Failure to derive a computational result is indicated by 'f' in the table.

| No. of displacement increments | Staggered | Monolithic | Hybrid $h_{\text{tol}} = 0.1$ | Hybrid $h_{\text{tol}} = 0.01$ |
| --- | --- | --- | --- | --- |
| 20  | f    | f | 296 | 94   |
| 100 | 1980 | f | 354 | 307  |
| 200 | 1667 | f | 463 | 547  |
| 500 | 1388 | f | 936 | 1203 |

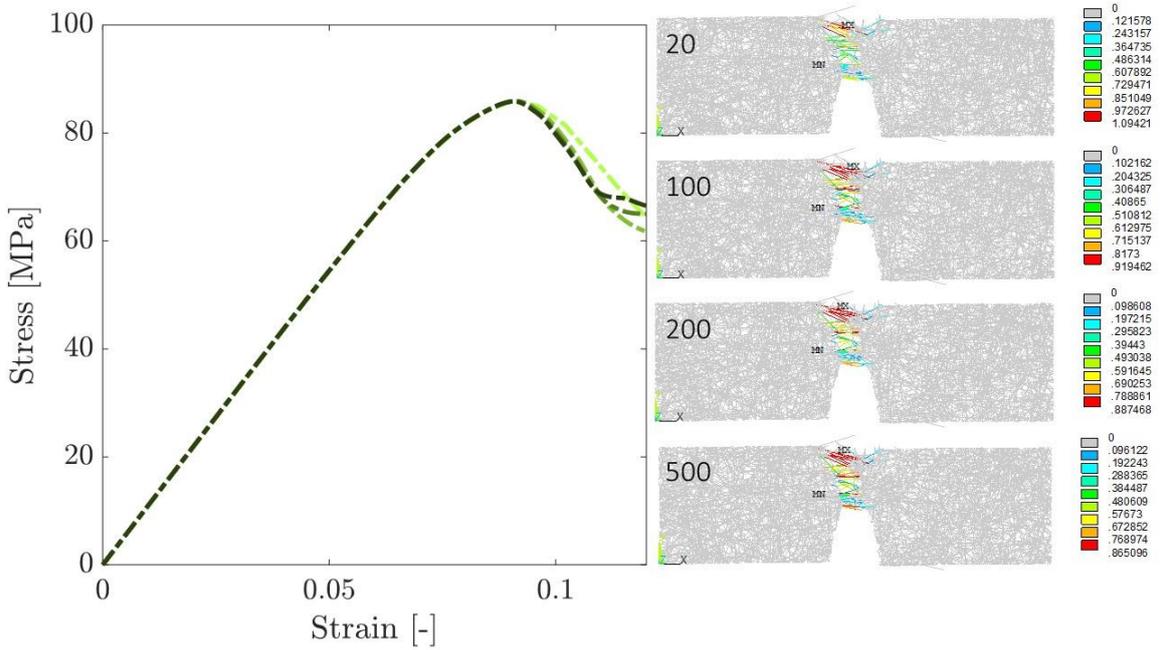

**Fig. 6**: Left: Stress-strain response of fiber networks with the density of $1000\ kg/m^3$ and fibers modeled as Timoshenko beams with the properties listed in **Table 3** and the facture energy $G_f = 0.1\ Nm$. The curves present results from the hybrid solution technique through the application of 20, 100, 200 and 500 (denoted by a light-to-dark color transition) displacement increments, respectively. Numerical effort is listed in **Table 4** and compared to the monolithic and staggered solution techniques. Right: Development of failure in the fibrous network as predicted through the prescription of 20, 100, 200 and 500 displacement increments, where the jump at complete rupture is $\xi_{max} = 0.85\ mm$. Any $\xi > \xi_{max}$ indicates the fiber has failed.

## 7    Summary and Conclusion

We studied the computational analysis of failure in fibrous materials, where the individual fibers are modeled as Timoshenko beams with embedded strong discontinuities. Representative benchmark examples have been used, where the recently proposed staggered solution method [39] has been tested against the monolithic solution strategy. Whilst the staggered approach is numerically robust, it does not use the consistent linearization of the nodal forces and therefore suffers from a poor convergence rate. This is especially the case for large displacement increments. The monolithic approach, in contrary, follows from the consistent linearization but can result in a non-positive definite element stiffness matrix, that then requires the solution of unstable equilibria. It is therefore practically not applicable to solve the benchmark problems studied in this work; it erratically fails, and step-size refinement is not always successful. We therefore proposed a novel hybrid solution method that forms the element stiffness through an adaptive 'mixing' of the stiffness of the monolithic and staggered approaches. It may also be seen as a matrix regularization technique to retain a positive definite element stiffness matrix while approaching the tangent stiffness matrix of the monolithic problem. The hybrid method results in a robust and computational efficient solution technique to explore failure in fibrous materials. The approach is general and may also accelerate the computational analysis of other failure problems.

**CRediT authorship contribution statement**

**Vedad Tojaga**: Conceptualization, Methodology, Software, Validation, Formal analysis, Investigation, Data curation, Writing – original draft, Writing – review & editing, Visualization, Project administration. **Artem Kulachenko**: Funding acquisition, Writing – review & editing, Supervision. **Sören Östlund**: Funding acquisition, Writing – review & editing, Supervision. **Thomas Christian Gasser**: Conceptualization, Methodology, Validation, Formal analysis, Investigation, Writing – original draft, Writing – review & editing, Visualization, Supervision.

**Acknowledgement**

V. Tojaga wishes to acknowledge Dr. August Brandberg for providing technical support in ANSYS and Professor Mijo Nikolic for insights into computational aspects.

This project has received funding from the European Union's Horizon 2020 research and innovation program under the Marie Skłodowska-Curie grant agreement No 764713, project FibreNet, and KTH Royal Institute of Technology.

The computations were performed on resources provided by the Swedish National Infrastructure for Computing (SNIC) at HPC2N (projects SNIC 2020/5-428 and SNIC 2021/6-51).